\documentclass[11pt]{article}
\usepackage{amsmath}
\usepackage{amssymb}  
\usepackage{epsfig}
\usepackage{amsthm}   
\usepackage[mathcal]{eucal}

\usepackage{url}

\newcommand{\F}{\mathcal{F}}
\newcommand{\Sf}{\mathcal{S}}
\newcommand{\J}{\mathcal{J}}

\newcommand{\R}{\mathbb{R}}
\newcommand{\BR}{\bar{\mathbb{R}}}

\DeclareMathOperator{\clconv}{cl\,conv}
\DeclareMathOperator{\ed}{dom}
\newcommand{\rot}{\Lambda}

\newcommand{\inner}[2]{\langle{#1},{#2}\rangle}

\newcommand{\norm}[1]{\|#1\|}
\newcommand{\normq}[1]{ {\|#1\|}^2 }

\newcommand{\tos}{\rightrightarrows} 

\newtheorem{theorem}{Theorem}[section]
\newtheorem{lemma}[theorem]{Lemma}
\newtheorem{corollary}[theorem]{Corollary}

\newtheorem{definition}[theorem]{Definition}

\title{On Gossez type (D) maximal monotone operators}

\author{M. Marques Alves\thanks{IMPA, Estrada Dona Castorina 110, 22460-320
    Rio de Janeiro, Brazil
   ({\tt maicon@impa.br})}\hspace{.5em}\thanks{Partially supported by Brazilian Ministry of Science and
    technology (MCT), scholarship PCI/DTI.}
  \and
    B. F. Svaiter\thanks{ IMPA, Estrada Dona Castorina 110, 22460-320 Rio de
    Janeiro, Brazil ({\tt benar@impa.br}) 
    tel: 55 (21) 25295112, fax: 55 (21)25124115.  }\hspace{.5em}
    \thanks{Partially supported by CNPq
    grants 300755/2005-8, 475647/2006-8 and by PRONEX-Optimization}
}

\date{}

\begin{document}

\maketitle

\begin{abstract}
  Gossez type (D) operators are defined in non-reflexive Banach spaces
  and share with the subdifferential a topological related property,
  characterized by bounded nets.
  In this work we present new properties and characterizations of
  these operators.
  The class (NI) was defined after Gossez defined the class (D) and
  seemed to generalize the class (D).
  One of our main results is the proof that these classes, type (D) and
  (NI), are identical.
  \\
\end{abstract}

\pagestyle{plain}


\section{Introduction}

Let $X$ be a real Banach space with topological dual $X^*$ and
topological bidual $X^{**}$. 
Whenever necessary, we will identify $X$ with its image under the
canonical injection of $X$ into $X^{**}$.
We denote by $\inner{\cdot}{\cdot}$ the
duality product in both $X\times X^*$ and $X^*\times X^{**}$,
\[
x^*(x)=\inner{x}{x^*},\qquad x^{**}(x^*)=\inner{x^*}{x^{**}}.
\]

A  point-to-set operator $T:X\tos X^*$ 
(respectively $T:X^{**}\tos X^{*}$)
 is a relation on $X$ to
$X^*$ 
(respectively on $X^{**}$ to $X^{*}$):
\[ 
 T\subset X\times X^* \quad\mbox{(respectively } T\subset
X^{**}\times X^*),
\]
and $r\in T(q)$ means $(q,r)\in T$.
An operator $T:X\tos X^*$, or $T:X^{**}\tos X^*$, is {\it monotone} if
\[
\inner{q-q'}{r-r'}\geq 0,\qquad \forall (q,r),(q',r')\in T.
\]
An operator $T:X\tos X^*$ is {\it maximal monotone} (in $X\times X^*$) if it
is monotone and maximal (whit respect to the inclusion) in the family
of monotone operators of $X$ into $X^*$.
An operator $T:X^{**}\tos X^*$ is maximal monotone (in $X^{**}\times
X^*$) if it is monotone and maximal (with respect to the inclusion) in
the family of monotone operators of $X^{**}$ into $X^*$.

Since the canonical injection of $X$ into $X^{**}$ allows one to
identity $X$ with a subset of $X^{**}$, any maximal monotone operator
$T:X\tos X^*$ is also a monotone operator from $X^{**}$ into $X^*$,
and admits one (or more) maximal monotone extension in $X^{**}$.
The problem of uniqueness of
such \emph{maximal monotone} extension to the bidual, has been
previously addressed by Gossez~\cite{gos-ope.jmaa71,gos-ran.pams72,
  gos-con.pams76,gos-ext.pams76}. Gossez defined the class of maximal
monotone operators of type (D), for which uniqueness of the extension
to the bidual is guaranteed~\cite{gos-ext.pams76}. This class of
maximal monotone operators has similar properties to the maximal
monotone operators in {\it reflexive} Banach spaces. In the reflexive
setting any maximal monotone operator is of Gossez type (D).
Properties like surjectivity of perturbations by duality mappings,
uniqueness of the extension to the bidual, etc., have been studied for
this class of operators by Gossez himself and also by others authors.

Attempting to generalize the Gossez type (D) class, Simons introduced
and studied in~\cite{sim-ran.jmaa96} the class of maximal monotone
operators of type (NI),
and proved that condition (NI) guarantees uniqueness of the extension
to the bidual and generalizes Gossez type (D) condition, but no
example showing that these classes are distinct was provided.
Equivalence between condition (D) and (NI) for \emph{linear}
(point-to-point) maximal monotone operators, bounded and unbouded was
proved in \cite{bau-bor-1999} and \cite{phe.sim-unb.jca98}
respectivelly.

Studding convex representations of maximal monotone operators in
non-reflexive spaces, the authors of the present paper
have found several properties of the class (NI) and new
characterizations in terms of the Fitzpatrick
family~\cite{Alves2008n,alv.sva-bro.jca08,alv.sva-max.jca09,alv.sva-new.jca09}.
In this work, making use of some of these tools, we will prove that a
maximal monotone operator is of type (NI) if and only if it is of
Gossez type (D).
Consequently, all properties studied recently in~\cite{Alves2008n,
  alv.sva-bro.jca08,alv.sva-max.jca09, alv.sva-new.jca09} can be
bridged to the class of operators of Gossez type (D).
We summarize them in Theorem~\ref{th:main}.
Beside that, 
 Lemma~\ref{lm:1}, 
and its consequence Theorem~\ref{th:a}, were important tools for us, but they
seem also to be relevants by themselves.
In Theorem~\ref{th:s} (and its corollary) we give a partial answer to
the question of equivalence between been of Gossez type (D) and having
a unique extension to the bidual.

\section{Basic definitions and classical results}
\label{sec:bd}
We use the notation $\norm{\cdot}$ for the norm in $X$, $X^{*}$ and $X^{**}$. 
In the cartesian product of Banach space we use the $\max$ norm.

The {\it effective domain} and the {\it epigraph} of $f:X\to \BR$ are defined,
respectively, as
\begin{align*}
 \ed(f)&=\{x\in X\,|\,f(x)<\infty\},\\
 \mbox{epi}(f)&=\{(x,\lambda)\in X\times \R\,|\,f(x)\leq \lambda\}.
\end{align*}
The function $f$ is \emph{convex} if its epigraph is convex, and it is proper
if $f>-\infty$ and $\ed(f)\neq\emptyset$.
For $\varepsilon\geq 0$, the $\varepsilon$-subdifferential of $f$ is
$\partial_\varepsilon f:X\tos X^*$,
\[
\partial_\varepsilon f(x)=\{ x^*\in X^*\;|\; f(y)\geq
f(x)+\inner{y-x}{x^*}-\varepsilon,\; \forall y\in X\}
\]  
and $\partial f$, the subdifferential of $f$, is the
$\varepsilon$-subdifferential for $\varepsilon=0$, that is $\partial f
=\partial_0 f$.

The $\varepsilon$-subdifferential of a proper convex
lower semicontinuous function is metricaly close to the subdiferential
of the function. This special property of the
$\varepsilon$-subdiferential, proved in~\cite{bro.roc-sub.pams65}, and called
\emph{Br\o ndsted-Rockafellar property} was a fundamental result used by
Rockafellar for proving the maximal monotonicity of the
subdifferential.
Br\o ndsted-Rockafellar property can be extended to maximal monotone
operators as follows:
\begin{definition}[\mbox{\cite{alv.sva-bro.jca08}}]
  A maximal monotone operator $T:X\tos X^*$ satisfies the \emph{strict
  Br\o ndsted-Rockafellar property} if, for any $\varepsilon>0$ and
  every $(x,x^*)$ such that
  \[ \inf_{(y,y^*)\in T} \inner{x-y}{x^*-y^*}> -\varepsilon 
  \]
  it holds that for any $\lambda>0$ there exists
  $(x_\lambda,x_\lambda^*)\in T$ such that
  \[ \norm{x-x_\lambda}<\lambda,\qquad
  \norm{x^*-x^*_\lambda}<\varepsilon/\lambda\;.
  \]
\end{definition}

Fenchel-Legendre conjugate of $f:X\to\BR$ is $f^*:X^*\to\BR$,
\[
f^*(x^*)=\sup_{x\in X}\inner{x}{x^*}-f(x).
\] 
Note that $f^*$ is always convex and lower semicontinuous in the
$\mbox{weak}^*$ (and hence strong) topology of $X^*$.
Under suitable assumptions, $f^{**}$ coincides with $f$ in $X$.
\begin{theorem}[Moreau]\label{th:mr}
  If $f:X\to\BR$ is a proper, convex, lower semicontinuous function, then
  $f^*$ is proper and $f^{**}(x)=f(x)$ for all $x\in X$.
\end{theorem}
An elegant proof of the above theorem is provided in~\cite{brezis}.
In the sequel we will need an auxiliary result, proved by Rockafellar inside the
proof of~\cite[Proposition 1]{roc-max.pjm70}.  For the sake of completeness, a proof
of this result is supplied in Appendix~\ref{ap:1}.

\begin{lemma}\label{lm:roc}
  Let $f:X\to\BR$ be a proper convex lower semicontinuous function.
  Then, for any $x^{**}\in X^{**}$ there exists $\alpha>0$ such
  that $f+\delta_\alpha$ is proper, lower semicontinuous and
 \[ f^{**}(x^{**})=(f+\delta_\alpha)^{**}(x^{**}),
 \]
 where $\delta_\alpha$ denotes the indicator function of
 $B_{X}[0,\alpha]$ in $X$,  $B_X[0,\alpha]$ being  the closed ball
   centered at $0$, with radius $\alpha$.
\end{lemma}

\section{Background materials}

In~\cite{gos-ope.jmaa71} Gossez introduced a class of maximal monotone operators in
\emph{non-reflexive} Banach spaces which share many properties with
maximal monotone operators defined in reflexive spaces.
This class, coined by Gossez as of type (D) operators, is the 
main concern of  this article.
\begin{definition}[\cite{gos-ope.jmaa71}]
Gossez's  monotone closure (with respect
to $X^{**}\times X^*$) of a maximal monotone operator $T:X\tos X^*$,
is $\overline T^{\,g}:X^{**}\tos X^*$,
\begin{equation}
  \label{eq:mc}
 \overline T^{\,g}=\{(x^{**},x^*)\in X^{**}\times X^* \,|\,
\inner{x^*-y^*}{x^{**}-y}\geq 0,\;\forall (y,y^*)\in T\}. 
\end{equation}
A maximal monotone operator $T:X\tos X^*$, is of Gossez type (D) if
for any $(x^{**},x^*)\in \overline T^{\,g}$, there exists a \emph{bounded}
net $\{(x_i,x_i^*)\}_{i\in I}$ in $T$ which 
converges to $(x^{**},x^*)$ in the $\sigma(X^{**},X^*)\times$strong
topology of $X^{**}\times X^*$.
\end{definition}
In~\cite{sim-ran.jmaa96} Simons introduced a what seemed to be a new
class of maximal monotone operators, called (NI), proved that this 
class encompasses Gossez type (D) class, and also generalized to this
class some previous results of Gossez for type (D) operators.
\begin{definition}[\cite{sim-ran.jmaa96}]
  \label{def:ni} A maximal monotone
  operator $T:X\tos X^*$ is of Simons type \emph{(NI)} if
  \[ \inf_{(y,y^{*})\in T}  \inner{y^*-x^*}{y-x^{**}}\leq 0,\qquad
   \forall (x^*,x^{**})\in X^*\times X^{**}.
   \]
\end{definition}
As we shall see, the classes Simons type (NI) and Gossez type (D) are identical.
We will also present new characterizations of the Gossez type (D)
operators.  Our main tool will be the Fitzpatrick functions which is
the next topic.

The Fitzpatrick function and the Fitzpatrick family of a maximal
monotone operator $T:X\tos X^*$ are defined, respectively,
as~\cite{fit-rep.pcmaanu88}:
\begin{equation}\label{fitz}
  \varphi_T:X\times X^*\to \BR,\qquad
  \varphi_T(x,x^*)
 =\sup_{(y,y^*)\in T}\,\inner{y}{x^*}+\inner{y^*}{x}-\inner{y}{y^*}
\end{equation}
and
\begin{equation} \label{eq:def.ft} \F_T=\left\{ h\in \BR^{X\times X^*}
    \left|
      \begin{array}{ll}
        h\mbox{ is convex and lower semicontinuous}\\
        \inner{x}{x^*}\leq h(x,x^*),\quad \forall (x,x^*)\in X\times X^*\\
        (x,x^*)\in T 
        \Rightarrow 
        h(x,x^*) = \inner{x}{x^*}
      \end{array}
    \right.
  \right\}.
\end{equation}
Next we summarize Fitzpatrick's results:
\begin{theorem}[\mbox{\cite[Theorem 3.10]{fit-rep.pcmaanu88}}] \label{th:fitz} 
  Let $T:X\tos X^*$ be maximal
  monotone. Then for any $h\in \F_T$  and $(x,x^*)\in X\times X^*$
  \[
  (x,x^*)\in T\iff h(x,x^*)=\inner{x}{x^*},
  \]
  and $\varphi_T$  is the smallest element
  of the  family $\F_T$.
\end{theorem}
Burachik and Svaiter defined and studied in~\cite{bur.sva-max.sva02}
the largest element of $\F_T$, the $\mathcal{S}$-function:
\begin{equation}\label{eq:def.s}
 \mathcal{S}_T=\sup_{h\in \F_T}\,h.
\end{equation}
Moreover, they also characterized it by
\begin{equation}\label{eq:char.s}
 \mathcal{S}_T=\clconv (p+\delta_T),
\end{equation}
where $\delta_T$ denotes the indicator function of $T$, $\clconv$ is
the convex lower semicontinuous closure in the strong topology of
$X\times X^*$ and $p:X\times X^*\to \R$ is the duality product,
\[
 p(x,x^*)=\inner{x}{x^*}.
\]
Burachik and Svaiter also proved that the family $\F_T$ is invariant under a
generalized conjugation, which we describe next.
Recall that $(X\times X^*)^*=X^*\times X^{**}$.
For $f:X\times X^*\to \BR$ \big($f:X^*\times X^{**}\to \BR$\big), define 
$\J f:X\times X^*\to \BR$ \big($\J_* f:X^*\times X^{**}\to \BR$\big) by
\[
 \J f(x,x^*)=f^*(x^*,x)\quad \big(\J_* f(x^*,x^{**})=f^*(x^{**},x^*)\big).
\]
To simplify the notation, define
 \begin{equation}
   \label{eq:def.rot}
 \rot:X^{**}\times X^*\to X^*\times X^{**},\quad
\rot(x^{**},x^*)=(x^*,x^{**}).  
 \end{equation}
Note that $\rot(X\times X^*)=X^*\times X$.
\begin{theorem}[\mbox{\cite[Theorem 5.3]{bur.sva-max.sva02}}]
 \label{th:x}
 Let $T:X\tos X^*\, \big(T:X^{**}\tos X^*\big)$ be maximal
 monotone. Then
 \[
  \J(\F_T)\subset \F_T \quad \big(\J_*(\F_{\rot T})\subset \F_{\rot T}\big)
 \]
\end{theorem}
In the sequel we will need a variant of~\cite[Theorem
3.4]{alv.sva-bro.jca08}, which we enunciate and prove next.
\begin{theorem}
  \label{th:variant}
  Let $h:X\times X^*\to \BR$ be a proper convex lower semicontinuous function.
Suppose that
\[
 h(x,x^*)\geq \inner{x}{x^*},
    \qquad \forall (x,x^*)\in X\times X^*\] 
\[ h^*(x^*,x^{**})\geq \inner{x^*}{x^{**}},
   \qquad \forall (x^*,x^{**})\in X^*\times X^{**}\] 
If $ h(z,z^*)-\inner{z}{z^*}\leq \varepsilon<\infty$ then, for any
$\eta,\mu>0$ such that $\eta\,\mu>1$, there exists $(y,y^*)\in
X\times X^*$ such that
\[ h(y,y^*)=\inner{y}{y^*},\quad
   \norm{y-z}\leq \eta \sqrt{\varepsilon}\;, \quad 
   \norm{y^*-z^*}\leq \mu\sqrt{\varepsilon}
\]
with strict inequalities if $\varepsilon>0$.
\end{theorem}
\begin{proof}
  If $\varepsilon=0$ this means that $h(z,z^*)=\inner{z}{z^*}$ and the pair
  $(y, y^*)=(z,z^*)$ has the desired properties.

  Suppose that $0<\varepsilon<\infty$. In this case, defining
  $  \tilde \varepsilon=\eta\mu\varepsilon$ we have
  \[
   h(z,z^*)-\inner{z}{z^*}<\tilde\varepsilon.
  \]
  Using now the above equation and Theorem 3.4 of \cite{alv.sva-bro.jca08},
  we conclude that for any $\lambda>0$ there
  exists $(y_\lambda,y_\lambda^*)$ such that
  \[   h(y_\lambda,y_\lambda^*)=\inner{y_\lambda}{y_\lambda^*},\quad
   \norm{y-z}< \lambda, \quad 
   \norm{y^*-z^*}< \tilde \varepsilon/\lambda.
   \]
   Now, the conclusion follows taking $\lambda=\eta\sqrt{\varepsilon}$.
\end{proof}

\section{Main results}

The next two results are the main tools of our analysis.
\begin{lemma}\label{lm:1}
  Take $f:X\times X^{*}\to\BR$ convex proper lower semicontinuous (in
  the strong topology) and define $\tilde f:X^{**}\times X^*\to\BR$,
  \[ 
  \tilde f(x^{**},x^*)=
  \begin{cases}
    f(x^{**},x^*),& x^{**}\in X,\\
    \infty &\mbox{otherwise}.
  \end{cases}
  \]
  Then, 
  \[ \clconv_{w-* \times s} \tilde
  f(x^{**},x^*)=f^{**}(x^{**},x^*),\qquad \forall (x^{**},x^*)\in
  X^{**}\times X^*
  \]
  where $\clconv_{w-*\times s}$ stands for the convex lower
  semicontinuous closure in the $\sigma(X^{**},X^*)\times$strong
  topology of $X^{**}\times X^*$.
  Equivalently,
  \[ f^{**}(x^{**},x^*)=\liminf_{(y,y^*)\to(x^{**},x^*)}f(y,y^*),
    \qquad \forall  (x^{**},x^*)\in
  X^{**}\times X^*,
  \]
  where the $\lim\inf$ is taken over all nets in $X\times X^*$ converging to
  $(x^{**},x^*)$ in the  $\sigma(X^{**},X^*)\times$strong
  topology of $X^{**}\times X^*$. 
\end{lemma}
\begin{proof}
  First use Theorem~\ref{th:mr} to conclude that $f^{**}$ and $\tilde
  f$ coincides on $X\times X^*$. Since $f^{**}(x^{**},x^*)\leq \tilde
  f(x^{**},x^*)$ for all $(x^{**},x^*)\in X^{**}\times X^*$ and
  $(x^{**},x^*)\mapsto f^{**}(x^{**},x^*)$ is lower semicontinuous in
  the $\sigma(X^{**},X^*)\times$ strong topology of $X^{**}\times
  X^*$, we have
 \[
 \clconv_{w-*\times s} \tilde f(x^{**},x^*)\geq
 f^{**}(x^{**},x^*),\qquad \forall (x^{**},x^*)\in X^{**}\times X^*.
 \]
 Suppose that the above inequality is strict at some $(\hat
 x^{**},\hat x^*)\in X^{**}\times X^*$.
 Define $g:X^{**}\times X^{*}\to \BR$ by
 $g(x^{**},x^*)=\clconv_{w-*\times s} \tilde f(x^{**},x^*)$.  Then,
 $\mbox{epi}(g)$ is convex, nonempty, and closed in the
 $\sigma(X^{**},X^*)\times$strong $\times$ strong topology of
 $X^{**}\times X^*\times \R$ and 
\[
\big((\hat x^{**},\hat x^*),f^{**}(\hat x^{**},\hat x^*)
\big)\notin  \mbox{epi}(g).
\]
Using the Hahn-Banach separation theorem in the linear space
$X^{**}\times X^*\times \R$, endowed with the locally convex topology
$\sigma(X^{**},X^*)\times$strong$\times$ strong, we have that there
exists a nonzero vector $(z^*,z^{**},\beta)\in X^*\times X^{**}\times
\R$ and $\mu\in \R$ such that
\begin{equation}
  \label{eq:0imp}
  \inner{x^{**}}{z^*}+\inner{x^*}{z^{**}}-\beta\lambda< \mu< \inner{\hat x^{**}}{z^*}+
 \inner{\hat x^*}{z^{**}}-\beta f^{**}(\hat x^{**},\hat x^*), 
\end{equation}
for all $(x^{**},x^*,\lambda)\in \mbox{epi}(g)$.  Since $\mbox{epi}(g)$
is nonempty and $g(x,x^*)\leq f(x,x^*)$ for all $(x,x^*)\in X\times
X^*$ we conclude, from the above inequality, that $\beta \geq
0$ and
\begin{equation}\label{eq:imp}
 \inner{x}{z^*}+\inner{x^*}{z^{**}}-\beta f(x,x^*)< \mu
\end{equation}
for all $(x,x^*)\in \ed f$.
Now, take $(\hat z^*,\hat z^{**})\in \ed(f^*)$ and
note that
\begin{equation}\label{eq:ca}
  \inner{x}{\hat z^*}+\inner{x^*}{\hat z^{**}}-f(x,x^*)\leq
  f^*(\hat z^*,\hat z^{**})<\infty
\end{equation}
for all $(x,x^*)\in X\times X^*$.
Multiplying the first inequality in the above equation by $\theta>0$
and adding the resulting inequality to the inequality
of~\eqref{eq:imp} we get
\[
\inner{x}{z^*+\theta \hat z^*}+\inner{x^*}{z^{**}+\theta \hat
  z^{**}}-(\beta+\theta)f(x,x^*)<\mu+\theta f^*(\hat z^*,\hat z^{**})
\]
for all $(x,x^*)\in \ed f$.
Dividing the above inequality by $\beta+\theta>0$, and defining
\[
z^*_\theta=\frac {1}{\beta+\theta}
(z^*+\theta\hat z^*),\;\;
z^{**}_\theta=\frac {1}{\beta+\theta}
(z^{**}+\theta\hat z^{**}),\;\;
\mu_\theta=\frac{\mu+\theta f^*(\hat z^*,\hat z^{**})}{\beta+\theta}
\]
we conclude that 
\[ f^*(z^*_\theta,z^{**}_\theta)\leq \mu_\theta.
\]
Therefore,
\[ 
f^{**}(\hat x^{**},\hat x^*)\geq \inner{z^*_\theta}{\hat
  x^{**}}+\inner{z^{**}_\theta}{\hat x^*} -\mu_\theta.
\]
Multiplication  of the above inequality by $\beta+\theta$ and direct
algebraic manipulation of the result yields
\[\inner{z^*+\theta\hat z^*}{\hat x^{**}}
+\inner{z^{**}+\theta\hat z^{**}}{\hat x^*}
\leq (\mu+\theta f^*(\hat z^*,\hat z^{**}))+(\beta+\theta)f^{**}(\hat
x^{**},\hat x^*).
\]
Taking the limit as $\theta\to 0^+$ on the above inequality we get
\[
\inner{z^*}{\hat x^{**}}+\inner{z^{**}}{\hat x^*}\leq 
\mu+\beta f^{**}(\hat x^{**},\hat x^*)
\]
which contradicts the second inequality of~\eqref{eq:0imp}.

The last part of the lemma follows from the first part and the
topological characterization of lower semicontinuous
closure by nets.
\end{proof}

\begin{theorem}\label{th:a}
  Let $f:X\times X^*\to\BR$ be a proper convex lower semicontinuous
  function. Then, for any $(x^{**},x^*)\in X^{**}\times X^{*}$ there
  exists a bounded net $\{(z_i,z^*_i)\}_{i\in I}$ in $X\times X^*$ which 
  converges to $(x^{**},x^*)$ in the
  $\sigma(X^{**},X^*)\times$strong topology of $X^{**}\times X^*$ and
   \[ f^{**}(x^{**},x^*)=\lim_{i\in I} f(z_i,z_i^*).\]
\end{theorem}
\begin{proof}
  Suppose that $f^{**}(x^{**},x^*)\in\R$.  First use Lemma~\ref{lm:roc} to conclude
  that there exists $\alpha>\norm{(x^{**},x^*)}$ such that
  $f+\delta_\alpha$ is proper, convex,
  lower semicontinuous and
\[  f^{**}(x^{**},x^*)=(f+\delta_\alpha)^{**}(x^{**},x^*).
\]
Therefore, using now Lemma~\ref{lm:1} we conclude that there exists
a net $\{(z_i,z^*_i)\}_{i\in I}$ in $X\times X^*$, converging to
$(x^{**},x^*)$ in the $\sigma(X^{**},X^*)\times$strong topology of
$X^{**}\times X^*$, such that
\[
\lim_{i\in I}(f+\delta_\alpha)(z_i,z^*_i)= f^{**}(x^{**},x^*).
\]
As the above limit is finite, the net can be chosen with
$(f+\delta_\alpha)(z_i,z^*_i)$ finite, which readily implies bondedness of
the net.

Suppose now that  $f^{**}(x^{**},x^*)=\infty$.
Existence of a bounded net $\{z_i\}_{i\in I}$ in $X$ converging to
$x^{**}$ in the $\sigma(X^{**},X^*)$ topology
follows from Goldstine theorem. 
Defining $z^*_i=x^*$ for all $i\in I$ we conclude,
using again Lemma~\ref{lm:1}, that $\liminf f(z_i,z_i^*)$ is
$\infty$. Hence, the net $\{(z_i,z_i^*)\}_{i\in I}$ has the desired properties.
\end{proof}
Recall that Fitzpatrick functions are defined in $X\times X^*$ and are
bounded bellow by the duality product. The next ressult, which is an
immediat consequence of Theorem~\ref{th:a}, deals with convex function
bouded bellow by the duality product.
\begin{corollary}
  \label{cr}
  If $f:X\times X^*\to\BR$ is convex lower semicontinuous,
  \[
   f(x,x^*)\geq \inner{x}{x^*}, \qquad \forall (x,x^*)\in X\times
  X^{*},
  \]
  and $\tilde f:X^{**}\times X^*\to \BR$ is defined as
  \[
  \tilde f(x^{**},x^*)=
  \begin{cases}
    f(x^{**},x^*),& x^{**}\in X,\\
    \infty &\mbox{otherwise},
  \end{cases}
  \]
  then
  \begin{equation}
    \label{eq:xxx}
        \clconv_{w-* \times s} \tilde
  f(x^{**},x^*)= f^{**}(x^{**},x^*)\geq \inner{x^*}{x^{**}}, \qquad \forall
   (x^*,x^{**})
   \in X^*\times
  X^{**}
 \end{equation}
  where $\clconv_{w-*\times s}$ stands for the convex lower
  semicontinuous closure in the $\sigma(X^{**},X^*)\times$strong
  topology of $X^{**}\times X^*$.
\end{corollary}
\begin{proof}
  The corollary hold trivially if $f\equiv \infty$. If $f$ is proper,
  apply Lemma~\ref{lm:1} to conclude that the equality
  in~\eqref{eq:xxx} holds. To prove the inequality in~\eqref{eq:xxx}
  us Theorem~\ref{th:a}.
\end{proof}

For $h:X\times X^*\to \BR$ and $(z,z^*)\in X\times X^*$ define
$h_{(z,z^*)}:X\times X^*\to \BR$ by~\cite{leg.sva-mon.sva05,alv.sva-new.jca09}
\[
 h_{(z,z^*)}(x,x^*)=h(x+z,x^*+z^*)-[\inner{x}{z^*}+\inner{z}{x^*}+\inner{z}{z^*}].
\]
By $R(T)$ we denote the {\it range} of $T:X\tos X^*$ and by
$(x,x^*)\in T(\cdot +z_0)$ we means $(x+z_0,x^*)\in T$. Moreover,
$\overline A$ denotes the closure (in the strong topology of $X^*$) of
a set $A\subset X^*$.

The subdifferential and the $\varepsilon$-subdifferential of the
function $\frac{1}{2}\|\cdot\|^2$ will be denoted by $J:X\tos X^*$ and
$J_\varepsilon:X\tos X^*$ respectively
\[
J(x)=\partial\; \frac{1}{2}\|x\|^2,\qquad J_\varepsilon
(x)=\partial_\varepsilon\; \frac{1}{2}\|x\|^2.
\]
The operator $J$ is known by the {\it duality mapping} of $X$. The operator
$J_{\varepsilon}$ was introduced by Gossez in~\cite{gos-ope.jmaa71}.

\begin{theorem}\label{th:main}
  Let $T:X\tos X^*$ be maximal monotone. Then all the conditions below
  are equivalent and they hold if and only if $T$ is of Gossez type
  (D):
  \begin{enumerate}
  \item $T$ is of type (NI);
  \item $(\Sf_T)^*(x^*,x^{**}) \geq\inner{x^*}{x^{**}},\qquad \forall
      (x^*,x^{**})\in X^*\times X^{**}$;
  \item For all $h\in\F_T$,
    \[ h^{*}(x^*,x^{**})\geq\inner{x^*}{x^{**}},\qquad \forall
      (x^*,x^{**})\in X^*\times X^{**};
    \]
  \item There exists $h\in\F_T$ such that
    \[ h^{*}(x^*,x^{**})\geq\inner{x^*}{x^{**}},\qquad \forall
      (x^*,x^{**})\in X^*\times X^{**};
    \]
  \item There exists  $h\in \F_T$ such that
    \[
    \inf_{(x,x^*)\in X\times X^*} h_{(x_0,x_0^*)}(x,x^*)+\frac{1}{2}\normq{x}+
    \frac{1}{2}
    \normq{x^*}=0,\qquad \forall (x_0,x_0^*)\in X\times X^*;
    \]
  \item  For \emph{all} $h\in \F_T$,
    \[
    \inf_{(x,x^*)\in X\times X^*} h_{(x_0,x_0^*)}(x,x^*)+\frac{1}{2}\normq{x}+
    \frac{1}{2}
    \normq{x^*}=0,\qquad \forall (x_0,x_0^*)\in X\times X^*;
    \]
 \item  $ \overline{R(T(\cdot+z_0)+ J)}=X^*$, for all 
    $z_0\in X$;
  \item $ \overline{R(T(\cdot+z_0)+ J_\varepsilon)}=X^*$, for all 
    $\varepsilon>0$, $z_0\in X$;
  \item $R(T(\cdot+z_0)+ J_\varepsilon)=X^*$, for all 
    $\varepsilon>0$, $z_0\in X$.
  \end{enumerate}
  Moreover, if any of these conditions hold, then
  \begin{description}
  \item[a)] $T$ satisfies the
  strict Br\o ndsted-Rockafellar property;
\item [b)] 
  $T$ has a unique maximal monotone extension to an operator from
  $X^{**}$ into $X^*$ and this extension is $\overline{T}^{\,g}:X^{**}\tos X^*$.
\item [c)] $(\Sf_T)^*=\varphi_{\Lambda\overline T^g}$\,;
\item  [d)] For any $h\in \F_T$, $ h^*\in \F_{\Lambda \overline T^g}$.
  \end{description}
\end{theorem}
\begin{proof}
  Equivalence between conditions 1-9 follows from~\cite[Proposition 1.3]{alv.sva-max.jca09},
~\cite[Theorem 1.2]{alv.sva-new.jca09}
  and~\cite[Theorem 3.6]{Alves2008n}.
Statements {\bf a)}, {\bf b)}, {\bf c)} and {\bf d)} have been proved in 
~\cite[Theorem 1.1]{alv.sva-max.jca09}.
The fact that a Gossez type (D) operator is of type (NI) has been proved
in~\cite[Lemma 15]{sim-ran.jmaa96}.

For finishing the proof, we will
prove that condition 1 implies that $T$ is of Gossez type (D).
For this aim, suppose $T$ is of type (NI) and 
\[
 (y^{**},y^*)\in\overline T^{\,g}.
\]
Take $h\in \F_T$. Using {\bf d)} we have $h^*\in \F_{\Lambda \overline
  T^{\,g}}$ and, since by Theorem~\ref{th:x}, $\J_*$ maps $\F_{\Lambda
  \overline T^{\,g}}$ into itself,
$\J_*(h^*)(y^*,y^{**})=\inner{y^*}{y^{**}}$. Therefore 
\begin{equation}\label{eq:e}
 h^{**}(y^{**},y^*)=\inner{y^*}{y^{**}}.
\end{equation}
By Theorem~\ref{th:a}, there exists a bounded net
$\{(z_i,z^*_i)\}_{i\in I}$ which converges to $(y^{**},y^*)$ 
in the $\sigma(X^{**},X^*)\times$strong  topology 
of $X^{**}\times X^*$ and
\begin{equation}\label{eq:h}
 h^{**}(y^{**},y^*)=\lim_{i\in I} h(z_i,z_i^*).
\end{equation}
Moreover, since $h^{**}(y^{**},y^*)$ is finite, the net can be chosen
so that $\{ h(z_i,z_i^*)\}_{i\in I}$ is bounded.
Let
\[ \varepsilon_i= h(z_i,z_i^*)-\inner{z_i}{z_i^*},\qquad i\in I.
\]
Bondedness of the net $\{(z_i,z^*_i)\}_{i\in I}$, together with its
convergence in the $\sigma(X^{**},X^*)\times$strong topology implies
that 
\[
\lim_{i\in I}\inner{z_i}{z_i^*}=\inner{y^{**}}{y^*}.
\]
Combining
this result with \eqref{eq:e}, \eqref{eq:h} 
we conclude that the real (non-negative) net $\{\varepsilon_i\}_{i\in
  I}$ is bounded and converges to $0$.

Since $h\in \F_T$ and $h^*\in\F_{\Lambda\overline T^g}$, the functions
$h$ and $h^*$ majorize the duality product in their respective
domains. Therefore, using the above definition of $\varepsilon_i$ and 
Theorem~\ref{th:variant} 
 for each $i\in I$, with $\eta_i=\sqrt{2}$
and $\mu_i=\sqrt{2}$,
we have that for each $i\in I$, there exists 
$(y_i,y_i^*)\in X\times X^*$ such that
\begin{equation}
  \label{eq:xxx}
  h(y_i,y_i^*)=\inner{y_i}{y_i^*},\quad\|z_i-y_i\|\leq \sqrt{ 2
  \varepsilon_i},\quad\|z_i^*-y_i^*\|\leq \sqrt{2 \varepsilon_i}.
\end{equation}
The inclusion $h\in \F_T$ and the firs equality in the above equation
readily implies
\[ (y_i,y_i^*)\in T,\qquad \forall i\in I.
\]
Using the boundedness of $\{(z_i,z^*_i)\}_{i\in I}$ and
$\{\varepsilon_i\}_{i\in I}$ together with the two inequalities
in~\eqref{eq:xxx} and we conclude that $\{(y_i, y^*_i)\}_{i\in I}$ is
bounded.
Moreover, since the nets $\{\norm{y_i-z_i}\}_{i\in I}$ and
$\{\norm{y^*_i-z^*_i}\}_{i\in I}$ converge to $0$, the net $\{(y_i,
y^*_i)\}_{i\in I}$ also converges to $(y^{**},y^*)$ in the
$\sigma(X^{**},X^*)\times$strong topology of $X^{**}\times X^*$.
\end{proof}

\begin{theorem}\label{th:s}
  If $T:X\tos X^*$ is maximal monotone and has a unique maximal
  monotone extension $\tilde T:X^{**}\tos X^*$ to the bidual, then
  one of the following conditions  holds:\\
  1) $T$ is of Gossez type (D);\\
  2) $T$ is affine and non-enlargeable, that is,
  \begin{equation}
    \label{eq:non.e}
     \varphi_T(x,x^*)=
  \begin{cases}
    \inner{x}{x^*},& (x,x^*)\in T,\\
    \infty&   \mbox{otherwise}.
  \end{cases}
  \end{equation}
 and $\F_T=\{\varphi_T\}$.
\end{theorem}

\begin{proof}
  Suppose $T$ is maximal monotone and has a unique maximal monotone
  extension $\tilde T:X^{**}\tos X^*$ to the bidual.
  Then from~\cite[Theorem 1.3]{alv.sva-max.jca09} we have either $T$
  is of type (NI), or $T$ is an affine manifold and
  $T=\ed(\varphi_T)$.
  Using Theorem~\ref{th:main} we
  obtain either condition 1 holds or  $T$ is an affine manifold and
  $T=\ed(\varphi_T)$.

  Now, suppose that $T$ is an affine manifold and
  $T=\ed(\varphi_T)$. Since $\varphi_T$ coincides with the duality
  product in $T$, we conclude that~\eqref{eq:non.e} holds. Let
  $h\in\F_T$. Since $\varphi_T\leq h$ and $h$ coincides with the duality
  product in $T$, we have $h=\varphi_T$.
\end{proof}

A direct consequence of the above theorem is that, for nonlinear
maximal monotone operators, being of Gossez type (D) is a necessary and
sufficient condition for uniqueness of the maximal monotone extension to
the bidual.
\begin{corollary}
  Suppose that $T:X\tos X^*$ is maximal monotone and $T$ is nonlinear,
  that is, $T$ is not an affine subspace of $X\times X^*$. Then $T$
  has a unique maximal monotone extension to the bidual if, and only
  if, $T$ is of Gossez type (D).
\end{corollary}
\section{Gossez contributions}

Since operators of type (NI) have been studied by many authors, it
seems relevant to stress some of Gossez's main contributions to the
study of type (D) operators.  This because the class of type (NI) is
equal to the class of type (D).

Gossez proved~\cite{gos-ope.jmaa71,gos-ran.pams72} that if $T:X\tos
X^*$ is of type (D):
\begin{enumerate}
\item $T$ has a unique maximal monotone extension to the bidual;
\item $R(\lambda J_\varepsilon+T)=X^*$, for all $\lambda,\varepsilon>0$;
\item $\overline{ R(T)}$ is convex;
\item if $D(T)=X$ and $T$ is coercive, then $R(T)$ is dense in
  $X^{*}$ in the $\sigma(X^*,X)$ topology.
\end{enumerate}
where 
\[ D(T)=\{x\in X\;|\; T(x)\neq\emptyset\}.
\]
Gossez also supplied an example~\cite{gos-ran.pams72} of a
\emph{linear} maximal monotone operator, non type (D), with a unique
maximal monotone extension to the bidual. This shows that
Theorem~\ref{th:s} can not be strengthened.

\appendix

\section{Proof of Lemma~\ref{lm:roc}}\label{ap:1}

For proving the Lemma~\ref{lm:roc} we will use the Fenchel-Rockafellar
duality formula:
\begin{theorem}\label{teo:fr}
  Let $f,g:X\to \BR$ be proper convex lower semicontinuous functions.
  Suppose that $f$ (or $g$) is continuous at a point $\hat x \in
  \ed(f)\cap \ed(g)$.  Then
\begin{equation} \label{eq:fr}
  \inf_{x \in X}\,f(x) + g(x)=\max_{x^* \in X^*}\,- f^*(x^*)-g^*(-x^*).
\end{equation}
\end{theorem}
For a proof of Theorem~\ref{teo:fr}, see~\cite{brezis}. This theorem
has been generalized by Attouch and Brezis~\cite{att-brez}, but we only need 
 its classical version.
\begin{proof}[proof of Lemma~\ref{lm:roc}]
For $\alpha>0$, define $\delta_\alpha:X\to \BR$ by
\[
 \delta_\alpha(x)=
 \begin{cases}
    0,& \|x\|\leq \alpha,\\
    \infty&   \mbox{otherwise}.
  \end{cases}
\]
Lower semicontinuity of $f+\delta_\alpha$ follows from the lower
semicontinuity of $f$ and $\delta_\alpha$.
Direct calculations yields, for any $x^*\in X^*$ and $\alpha>0$,
\begin{align}\label{eq:sl}
 \nonumber
 (f+\delta_\alpha)^*(x^*)&=\sup_{x\in X}\,\inner{x}{x^*}-(f+\delta_\alpha)(x)\\
 \nonumber               &=\sup_{x\in X}\,(x^*-\delta_\alpha)(x)-f(x)\\
                         &=-\inf_{x\in X}\,f(x)+(\delta_\alpha-x^*)(x).
 \end{align}
There exists $\hat x\in X$ such that $f(\hat x)\in\R$. Take 
\begin{equation}
  \label{eq:alpha}
  \alpha>\norm{\hat x}.
\end{equation}
In this case, $\hat x\in \ed(f)\cap \ed(\delta_\alpha)$
 and also $\delta_\alpha$ is continuous at $\hat x$.
 Now, let $x^*\in X^*$. Using~\eqref{eq:sl} and Theorem~\ref{teo:fr},
 for $f$ and $g=\delta_\alpha-x^*$ we get
\begin{align*}
 (f+\delta_\alpha)^*(x^*)&=-\inf_{x\in X}\,f(x)+(\delta_\alpha-x^*)(x)\\
                         &=-\max_{y^*\in X^*}\,-f^*(y^*)-(\delta_\alpha-x^*)^*(-y^*)\\
                         &=\min_{y^*\in X^*}\,f^*(y^*)+(\delta_\alpha)^*(x^*-y^*)\\
                         &= f^*\Box\, \delta_\alpha^*(x^*).
\end{align*}
Thus, $(f+\delta_\alpha)^*=f^*\Box\, \delta_\alpha^*$ and so that
\begin{equation}\label{eq:sl2}
 (f+\delta_\alpha)^{**}=(f^*\Box\, \delta_\alpha^*)^*=f^{**}+\delta_\alpha^{**},
\end{equation}
where the second identity follows from direct calculations.
Direct calculation yields, for $x^*\in X^*$, $x^{**}\in X^{**}$
\begin{equation}\label{eq:sl3}
\delta_\alpha^*(x^*)=\alpha\norm{x^*},\qquad
 \delta_\alpha^{**}(x^{**})=
 \begin{cases}
  0,& \|x^{**}\|\leq \alpha,\\
  \infty& \mbox{otherwise}.
 \end{cases}
\end{equation}
To end the proof, let $x^{**}\in X^{**}$. Taking
$\alpha>\max\{\norm{x^{**}},\norm{\hat x}\}$ and using~\eqref{eq:sl2}
and the expression in~\eqref{eq:sl3} for $\delta_\alpha^{**}$ we have
\[
 f^{**}(x^{**}) =
f^{**}(x^{**})+\delta_\alpha^{**}(x^{**})
=(f+\delta_\alpha)^{**}(x^{**}).
\]
\end{proof}


\end{document}